\documentclass[10pt]{article}

\usepackage{t1enc}
\usepackage{amsmath,amssymb,amsopn}
\usepackage{theorem}

\theorembodyfont{\sl}
\newtheorem{theorem}{Theorem}[section]
\newtheorem{definition}[theorem]{Definition}
\newtheorem{lemma}[theorem]{Lemma}

\theorembodyfont{\rm}\newtheorem{remark}[theorem]{Remark}
\theorembodyfont{\rm}

\begin{document}

\newcommand{\rot}{\mathop{\textrm{rot}}\nolimits}
\renewcommand{\div}{\mathop{\textrm{div}}\nolimits}
\newcommand{\grad}{\mathop{\textrm{grad}}\nolimits}
\renewcommand{\Re}{\mathop{\textrm{Re}}\nolimits}
\renewcommand{\Im}{\mathop{\textrm{Im}}\nolimits}
\newcommand{\rank}{\mathop{\textrm{rank}}\nolimits}

\title{On connections between domain specific constants in some norm inequalities}
\author{S{\'a}ndor Zsupp{\'a}n\\ \small{Berzsenyi D{\'a}niel Evang{\'e}likus Gimn{\'a}zium}\\
\small{H-9400, Sopron, Sz{\'e}chenyi t{\'e}r 11., Hungary}\\
\small{zsuppans@gmail.com}}
\maketitle

\begin{abstract}
We derive connections between optimal domain specific constants figuring in the Friedrichs-Velte inequality for conjugate harmonic functions, in the Babu\v{s}ka-Aziz inequality for the divergence and in the improved Poincar{\'e} inequality for the gradient.
With the same method we obtain for spatial domains an improved Poincar{\'e} inequality for the rotation in connection with the corresponding Babu\v{s}ka-Aziz inequality.
\end{abstract}



\section{Introduction}

In \cite{friedrichs1937} Friedrichs proved an inequality between the norms of square integrable conjugate harmonic functions on planar domains.
Horgan and Payne \cite{horganpayne1983} discovered that a smooth simply connected planar domain supports the Friedrichs inequality if and only if it supports the Babu\v{s}ka-Aziz inequality for the divergence \cite{babuska-aziz1972}, which ensures the stable solvability of the divergence equation in an appropriate function space on the domain.
Moreover, they proved an important equation involving the optimal domain specific constants figuring in the corresponding inequalities.
Velte \cite{velte1998} generalized this connection for smooth simply connected spatial domains and for the Babu\v{s}ka-Aziz inequality for the rotation using another variant of the Friedrichs inequality.
Costabel et al.\cite{costabeldauge2015} proved that this connection between the Friedrichs-Velte and Babu\v{s}ka-Aziz inequalities and constants remains valid without any smoothness assumptions on the domain and can be further generalized for differential forms on arbitrary dimensional domains, see \cite{costabel2015}.

These inequalities and constants are not only of theoretical but also of practical interest for the numerical solutions of problems in fluid dynamics and elasticity, see \cite{crouzeix1997,stoyan1999,stoyan2000} and references therein. 
Despite of their importance exact values of these constants are known in a few cases \cite{friedrichs1937, velte1998}, altough the inequalities are proved to be valid on general classes of domains.
Shapiro \cite{shapiro1980} proved the Friedrichs inequality for planar domains satisfying an interior cone condition, Acosta et al. \cite{acosta2006} established the Babu\v{s}ka-Aziz inequality for the divergence in the class of John domains, which is a generalization of the class of domains satisfying an interior cone condition.
Recently Jiang et al. \cite{jiang2014} proved that the validity of the Babu\v{s}ka-Aziz inequality is equivalent to the John condition and to the validity of an improved Poincar{\'e} inequality \cite{hurri1994} provided the domain satisfies the separation property (which condition is fulfilled for any simply connected plane domain).

In this paper we focus rather on the constants than on the inequalities itself.
Motivated by a result in \cite{duran2012} we establish a connection between the optimal constant in the improved Poincar{\'e} inequality and the Friedrichs-Velte constant, which connection implies the simultaneous validity of the corresponding inequalities.
With the same method using the Friedrichs constant connected to Babu\v{s}ka-Aziz constant for the rotation we derive an improved Poincar{\'e} inequality using the rotation instead of the gradient.

In section \ref{sec:prelim} we formulate the notation and the preliminaries.
Next in section \ref{sec:main} we derive the main result, Theorem \ref{thm:fv<->ip}, which states that domains satisfying the Hardy inequality simultaneously support the Friedrichs-Velte and the improved Poincar{\'e} inequalities.
As a byproduct we obtain upper estimations for the constant in the improved Poincar{\'e} inequality for star-shaped domains using known upper estimations for the Friedrichs-Velte constants from \cite{costabeldauge2015, horganpayne1983, payne2007}.
We also discuss geometric conditions for the problem domain in order to satisfy the conditions of Theorem \ref{thm:fv<->ip}.

\section{Notation and preliminary results}\label{sec:prelim}

In this paper $\Omega$ denotes a bounded domain in $\mathbb{R}^n$.
Let $L_2(\Omega)$ be the usual space of square integrable functions over $\Omega$.
The norm and the integral mean on $\Omega$ of $f\in L_2(\Omega)$ are $\|f\|^2=\int_{\Omega}|f|^2$ and $f_{\Omega}=\frac{1}{|\Omega|}\int_{\Omega}f$, respectively.
$H^1(\Omega)$ denotes the Sobolev space of functions with $\nabla f\in L_2(\Omega)^n$.
Its subspace $H_0^1(\Omega)$ is the closure of smooth functions with compact support in $\Omega$ under the norm $\|f\|_1^2=\|f\|^2+\|\nabla f\|^2$ of $H^1(\Omega)$.
$|f|_1=\|\nabla f\|$ denotes the seminorm on $H_0^1(\Omega)$ equivalent to the $H^1$-norm.

\begin{definition}\label{def:FVineq}
The domain $\Omega\subset\mathbb{R}^{n}$ ($n=2,3$) supports the Friedrichs-Velte inequality if there is a positive constant $\Gamma$ depending only on the domain $\Omega$ such that for every pair of square integrable conjugate harmonic functions $u$ and $v$ there holds
\begin{equation}\label{ineq:FV}
\|u\|^2\le \Gamma\|v\|^2\,\text{ provided }\,u_{\Omega}=0
\end{equation}
\end{definition}
The least possible constant in \eqref{ineq:FV} is denoted by $\Gamma_{\Omega}$ and is called the Friedrichs-Velte constants of the domain $\Omega$.
Conjugate harmonic means in Definition \ref{def:FVineq} the Cauchy-Riemann equations 
\begin{equation}\label{eq:CR}
\nabla u=\nabla^{\bot}v
\end{equation} 
with $\nabla^{\bot}=(\partial_2,-\partial_1)$ for planar domains and the Moisil-Teodorescu equations 
\begin{equation}\label{eq:MT}
\rot v=\nabla u \text{ and }\div v=0
\end{equation}
for spatial domains.
The normalization $u_{\Omega}=0$ means that $u$ belongs to the orthogonal complement of the kernel of the gradient in $L_2(\Omega)$.
The inequality \eqref{ineq:FV} was investigated by Friedrichs \cite{friedrichs1937} and Shapiro \cite{shapiro1980} for planar domains and then by Velte, \cite{velte1998} for simply connected spatial domains with sufficiently smooth boundary.
Velte \cite{velte1998} also formulated another related inequality for simply connected spatial domains estimating the norm of the vector valued function $v$ in \eqref{eq:MT} as
\begin{equation}
\label{ineq:FVt}\|v\|^2\le \tilde{\Gamma}_{\Omega}\|u\|^2
\text{ provided }\int_{\Omega}v\cdot\nabla\phi=0\text{ for every }\phi\in H^1(\Omega).
\end{equation}
This normalization means that the solution $v$ of \eqref{eq:MT} belongs to the orthogonal complement of the kernel of $\rot$ in $L_2(\Omega)^3$.
On simply connected domains it is equivalent to $v\cdot n$=0 on the boundary where $n$ denotes the unit normal vector.

The exact value of the Friedrichs-Velte constant does not depend on the size of $\Omega$ only on its shape.
Its value is known only for a few domains, \cite{costabeldauge2015, horganpayne1983, payne2007} contain useful estimations for the class of star-shaped domains.   
Costabel \cite{costabel2015} developed a generalization of the Friedrichs-Velte inequality for differential forms which generalization incorporates also the unification of \eqref{ineq:FV} and \eqref{ineq:FVt} along with \eqref{eq:CR} and \eqref{eq:MT} and it is valid at least in the class of Lipschitz domains.

As observed in \cite{costabeldauge2015, horganpayne1983, velte1998} the Friedrichs-Velte inequalities and constants are closely related to the Babu\v{s}ka-Aziz inequality and to the corresponding domain specific constants.
\begin{definition}\label{def:BAineqdiv}
The domain $\Omega\subset\mathbb{R}^n$ ($n=2,3$) supports the Babu\v{s}ka-Aziz inequality for the divergence if there is a positive constant $C$ depending only on the domain $\Omega$ such that for every $u\in L_{2}(\Omega)$ with $u_{\Omega}=0$ there is a $v\in H_0^1(\Omega)^n$ such that $\div v=u$ and
\begin{equation}\label{ineq:BAdiv}
|v|_{1}^2\le C\|u\|^2.
\end{equation}
\end{definition}
The least possible constant in \eqref{ineq:BAdiv} is denoted by $C_{\Omega}$ and is called the Babu\v{s}ka-Aziz constant for the divergence of the domain.
$C_{\Omega}<\infty$ was proved for bounded Lipschitz domains in \cite{babuska-aziz1972} and was generalized for John domains in \cite{acosta2006}.

The inequality \eqref{ineq:BAdiv} can be formulated as $|v|_{1}^2\le C\|\div v\|^2$ for every function $v$ in the orthogonal complement of the kernel of the divergence in $H_0^1(\Omega)^n$.
Similarly there is a Babu\v{s}ka-Aziz inequality for the rotation:
\begin{equation}\label{ineq:BArot}
|v|_{1}^2\le \tilde{C}_{\Omega}\|\rot v\|^2
\end{equation}
provided $v$ is in the orthogonal complement of the kernel of the rotation in $H_0^1(\Omega)^3$.
According to \cite{costabeldauge2015,costabel2015, velte1998} there also hold
\begin{equation}\label{eq:BA-FV}
C_{\Omega}=1+\Gamma_{\Omega}\text{ and }\tilde{C}_{\Omega}=1+\tilde{\Gamma}_{\Omega}
\end{equation}
for any planar or spatial domain the constants being simultaneously finite or infinite.

The third class of inequalities utilized in this paper is the class of the improved Poincar{\'e} inequalities.

\begin{definition}\label{def:iPineq}
The domain $\Omega\subset\mathbb{R}^n$ supports the improved Poincar{\'e} inequality if there is a positive constant $P$ depending only on the domain $\Omega$ and on the exponents $\alpha$, $p$ and $q$ such that
\begin{equation}\label{ineq:pq_iP}
\|u-u_{\Omega}\|_{L_q(\Omega)}^q\le P\|d_{\Omega}^{\alpha}\nabla u\|_{L_p(\Omega)}^q
\end{equation}
holds for every $u\in L_{1,loc}(\Omega)$ such that $d^{\alpha}\nabla u\in L_p(\Omega)$, where $d_{\Omega}(x)=\text{dist}(x,\partial\Omega)$ is the distance of $x\in\Omega$ to the boundary.
\end{definition}
It was proved in \cite{boas1988} in case $p=q$ for domains whose boundary is locally the graph of a H{\"o}lder continuous function of order $\alpha$ and it was generalized in \cite{hurri1994} for $0\le\alpha\le 1$, $p(1-\alpha)<n$ and $p\le q\le\frac{np}{n-p(1-\alpha)}$ in a class of domains including John-domains.

In case $\alpha=0$ and $q=p=2$ one has the classical Poincar{\'e} inequality but in this paper we utilize the case $q=p=2$ and $\alpha=1$ in the form
\begin{equation}\label{ineq:iP}
\|u-u_{\Omega}\|^2\le P_{\Omega}\|d_{\Omega}\nabla u\|^2,
\end{equation}
wherein the improved Poincar{\'e} constant $P_{\Omega}$ of $\Omega$ is the least possible positive constant satisfying \eqref{ineq:iP}.
For bounded simply connected planar domains it was proved in \cite{jiang2014} that the domain supports the Babu\v{s}ka-Aziz inequality \eqref{ineq:BAdiv} iff it supports the improved Poincar{\'e} inequality \eqref{ineq:iP} and iff $\Omega$ is a John domain. For more general domains there are additional properties needed in order to have equivalence between $\Omega$ being a John domain and $\Omega$ supporting the inequalities \eqref{ineq:BAdiv} and \eqref{ineq:iP}, c.f. \cite{jiang2014}.

\begin{remark}
As proved in \cite{chuawheeden2010} for convex domains in arbitrary dimensions the constant $P$ in \eqref{ineq:pq_iP} can be estimated from above by a scalar multiple of the product $\eta(\Omega)^{2\alpha}\text{diam}(\Omega)^{2-2\alpha}$, where $\eta(\Omega)$ and $\text{diam}(\Omega)$ denote the eccentricity and the diameter of $\Omega$, respectively.
This estimator is independent of $\text{diam}(\Omega)$ only for $\alpha=1$ in which case the improved Poincar{\'e} constant $P_{\Omega}$ in \eqref{ineq:iP} can be estimated by a scalar multiple (depending only on the dimension $n$) of  $\eta(\Omega)^{2}$.
Hence if one wants to derive a correspondence between the diameter invariant Friedrichs-Velte constant and the improved Poincar{\'e} constants then one has to choose $\alpha=1$ in \eqref{ineq:pq_iP}. 
\hfill$\Box$
\end{remark}

\section{Connection between the constants}\label{sec:main}

In this section, motivated by Theorem 5.1 an Theorem 5.3 in \cite{duran2012}, we derive connections between the Friedrichs-Velte and improved Poincar{\'e} inequalities and constants.

\begin{lemma}\label{lem:ip->fv}
If the domain $\Omega\subset\mathbb{R}^n$ ($n=2,3$) supports the improved Poincar{\'e} inequality \eqref{ineq:iP} with the constant $P_{\Omega}$, then it also supports the Friedrichs-Velte inequality \eqref{ineq:FV} with the constant
\begin{equation}\label{ineq:iP->FV}
\Gamma_{\Omega}\le 4P_{\Omega}.
\end{equation} 
\end{lemma}
\textbf{Proof.}
Let $u$ and $v$ be conjugate harmonic functions in the sense of \eqref{eq:MT} on the spatial domain $\Omega$ supporting the improved Poincar{\'e} inequality.
We develop the norm on the right-hand side of \eqref{ineq:iP}:
\begin{eqnarray*}
\|d_{\Omega}\nabla u\|^2&=&\int_{\Omega}d_{\Omega}^2|\nabla u|^2=
\int_{\Omega}d_{\Omega}^2\nabla u\cdot\rot v
=\int_{\Omega}v\cdot\rot\left(d_{\Omega}^2\nabla u\right)\\
&=&\int_{\Omega}v\cdot\left(2d_{\Omega}\nabla d_{\Omega}\times\nabla u\right)
=\int_{\Omega}2d_{\Omega}\nabla u\cdot\left(v\times\nabla d_{\Omega}\right)
\end{eqnarray*}
We estimate by the Cauchy-Schwarz inequality and we also use $\left|\nabla d_{\Omega}\right|=1$ for the boundary distance function of $\Omega$ valid a.e. in $\Omega$.
\begin{equation}
\left|\int_{\Omega}2d_{\Omega}\nabla u\cdot\left(v\times\nabla d_{\Omega}\right)\right|
\le 2\left(\int_{\Omega}d_{\Omega}^2|\nabla u|^2\right)^{\frac{1}{2}}\left(\int_{\Omega}|v|^2\right)^{\frac{1}{2}}
\end{equation}
There follows
\begin{equation}
\|d_{\Omega}\nabla u\|^2\le 2\|v\|\cdot\|d_{\Omega}\nabla u\|.
\end{equation}
By the improved Poincar{\'e} inequality there follows
\begin{equation}
\|u-u_{\Omega}\|^2\le P_{\Omega}\|d_{\Omega}\nabla u\|^2
\le 4P_{\Omega}\|v\|^2,
\end{equation}
which implies \eqref{ineq:iP->FV} for the Friedrichs-Velte constant $\Gamma_{\Omega}$ of the domain.

For planar domains the above proof remains valid with minor changes.
Let $u$ and $v$ be conjugate harmonic functions on the planar domain $\Omega$ in the sense of the Cauchy-Riemann equations \eqref{eq:CR}. 
\begin{eqnarray*}
\|d_{\Omega}\nabla u\|^2&=&\int_{\Omega}d_{\Omega}^2|\nabla u|^2=
\int_{\Omega}d_{\Omega}^2\nabla u\cdot\nabla^{\bot} v
=\int_{\Omega}v\cdot\rot\left(d_{\Omega}^2\nabla u\right)\\
&=&-\int_{\Omega}v\left(2d_{\Omega}\nabla^{\bot}d_{\Omega}\cdot\nabla u\right)
=-\int_{\Omega}\left(2d_{\Omega}\nabla u\right)\cdot\left(v\nabla^{\bot}d_{\Omega}\right)
\end{eqnarray*}
We use $\left|\nabla^{\bot}d_{\Omega}\right|=|\nabla d_{\Omega}|=1$ valid a.e. in $\Omega$ and estimate by the Cauchy-Schwarz inequality.
\begin{equation*}
\left|\int_{\Omega}\left(2d_{\Omega}\nabla u\right)\cdot\left(v\nabla^{\bot}d_{\Omega}\right)\right|
\le 2\left(\int_{\Omega}|v|^2\right)^{\frac{1}{2}}\left(\int_{\Omega}d_{\Omega}^2|\nabla u|^2\right)^{\frac{1}{2}}
\end{equation*}
This implies the same inequality as in the spatial case: $\Gamma_{\Omega}\le 4P_{\Omega}$.
\hfill$\Box$

\begin{remark}\label{rem:duran}
Theorem 5.1 in \cite{duran2012} formulates a similar result between the Babu\v{s}ka-Aziz constant for the divergence and the improved Poincar{\'e} constant for arbitrary dimensional domains, however, with another dimension dependent constant which value is not specified. This reads with the notation of Lemma \ref{lem:ip->fv} 
\begin{equation}\label{ineq:duran}
\sqrt{1+\Gamma_{\Omega}}\le c_n\cdot\left(1+\sqrt{P_{\Omega}}\right).
\end{equation}
Altough Lemma \ref{lem:ip->fv} is proved only for planar and spatial domains but we now have an explicit constant in the inequality estimating $\Gamma_{\Omega}$ by $P_{\Omega}$ from above.
\hfill$\Box$
\end{remark}

\begin{remark}\label{rem:ip->fv_distance}
In the proof of Lemma \ref{lem:ip->fv} we used $\left|\nabla d_{\Omega}\right|=1$ for the boundary distance function valid a.e. in $\Omega$. The proof remains valid using instead of $d_{\Omega}$ another weight function $w$ with bounded gradient on $\Omega$, $|\nabla w|\le c$ and with zero boundary values on $\partial\Omega$ (appropiate solution of the eikonal equation). 
However, in this case one has to use another weighted Poincar{\'e} inequality instead of \eqref{def:iPineq} and one has another constant in \eqref{ineq:iP->FV} instead of 4.
\hfill$\Box$
\end{remark}

\begin{remark}\label{rem:polygons}
The constant 4 in \eqref{ineq:iP->FV} can be improved, for example for convex polygons one has more: $\Gamma_{\Omega}\le P_{\Omega}$.
In order to prove this we consider first that by \cite{gustafsson1998} every convex polygon has a unique mother body (skeleton) consisting of line segments which are subsets of bisectors of angles between appropriate two sides of the polygon.
The convex polygon $\Omega$ has a partition along this mother body $\Omega=\cup_j\Omega_j$.
We have $|\nabla d_{\Omega}|=1$ and $\Delta d_{\Omega}=0$ in each subpolygon $\Omega_j$ hence there follows
\begin{equation*}
\Delta\left(d_{\Omega}^{2}\right)=2|\nabla d_{\Omega}|^2+2d_{\Omega}\Delta d_{\Omega}=2
\end{equation*}
for the Laplacian of the square of the boundary distance function of $\Omega$ in each subpolygon $\Omega_j$.
On the boundary of each subpolygon $\Omega_j$ we have $\frac{\partial d_{\Omega}}{\partial n_j}>0$ because $\Omega$ is convex.
By partial integration we have on each $\Omega_j$
\begin{eqnarray*}
\int_{\Omega_j}d_{\Omega}^2|\nabla u|^2&=&\int_{\Omega_j}d_{\Omega}^2|\nabla^{\bot} v|^2=\int_{\Omega_j}d_{\Omega}^2|\nabla v|^2=\int_{\Omega_j}d_{\Omega}^2\Delta\left(\frac{1}{2}v^2\right)=\\
&=&\int_{\partial\Omega_j}d_{\Omega}v\left(d_{\Omega}\frac{\partial v}{\partial n_j}-v\frac{\partial d_{\Omega}}{\partial n_j}\right)+\frac{1}{2}\int_{\Omega_j}v^2\Delta\left(d_{\Omega}^{2}\right)
\end{eqnarray*}
for a conjugate harmonic pair $u$ and $v$.
Now summing up all these equations some boundary terms cancel out and we obtain
\begin{equation*}
\|d_{\Omega}\nabla u\|^2=\int_{\Omega}v^2-\sum_j\int_{\partial\Omega_j}d_{\Omega}v^2\frac{\partial d_{\Omega}}{\partial n_j},
\end{equation*}
Using $\frac{\partial d_{\Omega}}{\partial n}>0$ on the boundaries of the subpolygons, there follows
\begin{equation*}
\|d_{\Omega}\nabla u\|^2\le\|v\|^2,
\end{equation*}
which gives $\Gamma_{\Omega}\le P_{\Omega}$.
\hfill$\Box$
\end{remark}

Lemma \ref{lem:ip->fv} says that each domain which supports the improved Poincar{\'e} inequality also supports the Friedrichs-Velte inequality.
In order to prove the reverse direction we use 
\begin{itemize}
\item the equality $C_{\Omega}=1+\Gamma_{\Omega}$ proved in \cite{costabeldauge2015,costabel2015} stating the simultaneous finiteness of the Friedrichs-Velte and the Babu\v{s}ka-Aziz constants for a domain $\Omega$ without assuming any boundary regularity and

\item Theorem 5.3 in \cite{duran2012}, in which the finiteness of the improved Poincar{\'e} constant was proved assuming the finiteness of the Babu\v{s}ka-Aziz constant $C_{\Omega}$ provided the domain $\Omega$ also supports a Hardy type inequality involving the boundary distance function.
\end{itemize}

\begin{definition}\label{def:Hardyineq}
The bounded domain $\Omega\subset\mathbb{R}^n$ supports the Hardy inequality if there is a finite positive constant H depending only on $\Omega$ such that
\begin{equation}\label{ineq:hardy}
\int_{\Omega}\frac{u^2}{d_{\Omega}^2}\le H\int_{\Omega}|\nabla u|^2
\end{equation}
holds for every $u\in H_0^1(\Omega)^n$.
The Hardy constant $H_{\Omega}$ of the domain is the least positive constant $H$ for which \eqref{ineq:hardy} holds.
\end{definition}

\begin{lemma}\label{lem:fv->ip}
If the domain $\Omega\subset\mathbb{R}^n$ ($n=2,3$) supports the Friedrichs-Velte inequality \eqref{ineq:FV} with the constant $\Gamma_{\Omega}$ and the Hardy inequality \eqref{ineq:hardy} with the constant $H_{\Omega}$, then it also supports the improved Poincar{\'e} inequality \eqref{def:iPineq} with the constant
\begin{equation}\label{ineq:FV->iP}
P_{\Omega}\le H_{\Omega}(1+\Gamma_{\Omega}).
\end{equation}
\end{lemma}
\textbf{Proof.}
The proof is due to Dur{\'a}n \cite{duran2012}, Theorem 5.3, which we reproduce here for the convenience
of the reader using the notation of the present paper.
Given $u\in H^1(\Omega)$ with zero integral mean $u_{\Omega}=0$ let $v\in H_0^1(\Omega)^n$ such that
\begin{equation}
\div v=u\text{ and }\|\nabla v\|^2\le C_{\Omega}\|u\|^2.
\end{equation}
Using the Hardy inequality \eqref{ineq:hardy} there follows
\begin{equation}
\|u\|^2=\int_{\Omega}u\div v=-\int_{\Omega}v\cdot\nabla u\le\left\|\frac{v}{d_{\Omega}}\right\|\cdot\left\|d_{\Omega}\nabla u\right\|\le H_{\Omega}^{\frac{1}{2}}\|\nabla v\|\cdot\left\|d_{\Omega}\nabla u\right\|.
\end{equation}
Using now the Babu\v{s}ka-Aziz inequality there follows
\begin{equation}
\|u\|^2\le H_{\Omega}^{\frac{1}{2}}C_{\Omega}^{\frac{1}{2}}\|u\|\cdot\left\|d_{\Omega}\nabla u\right\|
\end{equation}
which gives \eqref{ineq:FV->iP} utilizing $C_{\Omega}=1+\Gamma_{\Omega}$.
\hfill$\Box$

\begin{remark}\label{rem:hardynecessity}
Lemma \ref{lem:fv->ip} means that $H_{\Omega}<\infty$ is sufficient for an estimation of the improved Poincar{\'e} constant by the Friedrichs-Velte constant.
Example 4.1 in \cite{jiang2014} shows that the validity of the Hardy inequality \eqref{ineq:hardy} is not necessary for such an estimation.
\hfill$\Box$
\end{remark}

The proofs of Lemma \ref{lem:ip->fv} and Lemma \ref{lem:fv->ip} are applicable with minor changes for the other Friedrichs-Velte inequality \eqref{ineq:FVt} and for the Babu\v{s}ka-Aziz inequality \eqref{ineq:BArot} for the rotation.

\begin{lemma}\label{lem:fvt->iprot}
If the domain $\Omega\subset\mathbb{R}^3$ supports the Hardy inequality \eqref{ineq:hardy} and the Friedrichs-Velte inequality \eqref{ineq:FVt} with the finite constants $H_{\Omega}$ and $\tilde{\Gamma}_{\Omega}$,respectively, then there is a finite positive constant $\tilde{P}$ depending only on $\Omega$ such that the inequality
\begin{equation}\label{ineq:iProt}
\|v\|^2\le \tilde{P}\|d_{\Omega}\rot v\|^2
\end{equation}
is valid for every $v$ in the orthogonal complement of the kernel of $\rot$ in $L_2(\Omega)^3$.
Moreover, we have $\tilde{P}_{\Omega}\le H_{\Omega}(1+\tilde{\Gamma}_{\Omega})$ for the least possible constant $\tilde{P}_{\Omega}$ in \eqref{ineq:iProt}.
\end{lemma}
\textbf{Proof.}
The proof is essentially the same as that of Lemma \ref{lem:fv->ip}.
First use $\tilde{C}_{\Omega}=1+\tilde{\Gamma}_{\Omega}$, see \cite{costabel2015,velte1998}.
By the Babu\v{s}ka-Aziz inequality for the rotation there exists $w\in H_0^1(\Omega)^3$ such that
\begin{equation}
v=\rot w\text{ and }\|\nabla w\|^2\le \tilde{C}_{\Omega}\|v\|^2.
\end{equation}
There follows by the Hardy inequality \eqref{ineq:hardy}
\begin{equation}
\|v\|^2=\int_{\Omega}v\cdot\rot w=\int_{\Omega}w\cdot\rot v\le\left\|\frac{w}{d_{\Omega}}\right\|\cdot\left\|d_{\Omega}\rot v\right\|\le H_{\Omega}^{\frac{1}{2}}\|\nabla w\|\cdot\left\|d_{\Omega}\rot v\right\|.
\end{equation}
Substituting now the Babu\v{s}ka-Aziz inequality for the rotation gives
\begin{equation}
\|v\|^2\le H_{\Omega}^{\frac{1}{2}}\tilde{C}_{\Omega}^{\frac{1}{2}}\|v\|\cdot\left\|d_{\Omega}\rot v\right\|,
\end{equation}
from which the statement of the lemma follows.
\hfill$\Box$

\begin{remark}\label{rem:poi_for_rot}
In the improved Poincar{\'e} inequality \eqref{ineq:iP} the function $u-u_{\Omega}$ belongs to the orthogonal complement of the kernel of the gradient, hence \eqref{ineq:iP} can be also stated as $\|u\|^2\le P_{\Omega}\|d_{\Omega}\nabla u\|^2$ for every $u\in\left(\mathrm{ker}\nabla\right)^{\bot}$.
The previous Lemma \ref{lem:fvt->iprot} formulates an analogous inequality for the rotation instead of the gradient, hence \eqref{ineq:iProt} can be seen as an improved Poincar{\'e} inequality for the rotation.
According to \cite{costabel2015} we have $\tilde{\Gamma}_{\Omega}<\infty$ at least for Lipschitz domains, for which the Hardy constant $H_{\Omega}$ is also finite \cite{necas1962}, hence \eqref{ineq:iProt} is valid for spatial Lipschitz domains.
\hfill$\Box$
\end{remark}

\begin{remark}\label{rem:2Dcase}
In two dimensions we have $\|d_{\Omega}\nabla u\|=\|d_{\Omega}\nabla^{\bot} u\|$, hence the improved Poincar{\'e} constant for the vector-curl $\nabla^{\bot}$ (adjoint of the scalar rotation) coincides with the usual Poincar{\'e} constant for the gradient.
For planar domains the Friedrichs-Velte constants $\Gamma_{\Omega}$ and $\tilde{\Gamma}_{\Omega}$ are equal as well.
\hfill$\Box$
\end{remark}

The counterpart of Lemma \ref{lem:ip->fv} is the following
\begin{lemma}\label{lem:iprot->fvt}
If the domain $\Omega\subset\mathbb{R}^3$ supports the inequality \eqref{ineq:iProt} for every $v$ in the orthogonal complement of the kernel of the $\rot$ in $L_2(\Omega)^3$ with the least possible positive constant $\tilde{P}_{\Omega}$, then it also supports the Friedrichs-Velte inequality \eqref{ineq:FVt}.
Moreover, we have $\tilde{\Gamma}_{\Omega}\le 4\tilde{P}_{\Omega}$.
\end{lemma}
\textbf{Proof.}
The proof is similar to that of Lemma \ref{lem:ip->fv}.
Let $u$ and $v$ be conjugate harmonic functions on the spatial domain $\Omega$ in the sense of the Moisil-Teodorescu equations \eqref{eq:MT} normalized such that $v$ lies in the orthogonal complement of the kernel of $\rot$ in $L_2(\Omega)^3$.
We develop the norm on the right-hand side of \eqref{ineq:iProt} using $|\nabla d_{\Omega}|=1$ a.e. in $\Omega$.
\begin{eqnarray*}
\|d_{\Omega}\rot u\|^2&=&\int_{\Omega}d_{\Omega}^2|\rot u|^2=
\int_{\Omega}d_{\Omega}^2\nabla u\cdot\rot v
=\int_{\Omega}u\div\left(d_{\Omega}^2\rot v\right)\\
&=&\int_{\Omega}u\left(2d_{\Omega}\nabla d_{\Omega}\cdot\rot v\right)
\le 2\|u\|\cdot\|d_{\Omega}\rot u\|
\end{eqnarray*}
This implies by \eqref{ineq:iProt}
\begin{equation}
\|v\|^2\le 4\tilde{P}_{\Omega}\|u\|^2,
\end{equation}
which means the Friedrichs-Velte inequality \eqref{ineq:FVt} with $\tilde{\Gamma}_{\Omega}\le 4\tilde{P}_{\Omega}$. 
\hfill$\Box$

\begin{remark}\label{rem:unification}
Considering the similarities between the proofs of Lemma \ref{lem:fv->ip} and Lemma \ref{lem:fvt->iprot} and between their reversed counterparts Lemma \ref{lem:ip->fv} and Lemma \ref{lem:iprot->fvt} they could be possibly unified in the framework of \cite{costabel2015} provided there is a usable version of the utilized Hardy and improved Poincar{\'e} inequalities for differential forms.  However, the investigation of this is beyond the scope of the present paper.
\hfill$\Box$
\end{remark}

As a consequence of these lemmata we obtain the following
\begin{theorem}\label{thm:fv<->ip}
If the bounded domain $\Omega$ supports the Hardy inequality \eqref{ineq:hardy}, then $\Omega$ supports the Friedrichs-Velte \eqref{ineq:FV} and \eqref{ineq:FVt} and simultaneously the Babu\v{s}ka-Aziz inequalities \eqref{ineq:BAdiv} and \eqref{ineq:BArot} if and only if $\Omega$ supports the improved Poincar{\'e} inequalities \eqref{ineq:iP} and \eqref{ineq:iProt}, respectively.
Moreover, we have 
\begin{equation}\label{ineq:fviphconstants}
\frac{1}{4}\Gamma_{\Omega}\le P_{\Omega}\le H_{\Omega}\left(1+\Gamma_{\Omega}\right)
\text{ and }
\frac{1}{4}\tilde{\Gamma}_{\Omega}\le \tilde{P}_{\Omega}\le H_{\Omega}\left(1+\tilde{\Gamma}_{\Omega}\right)
\end{equation}
for the domain spacific constants in the corresponding inequalities.
\hfill$\Box$
\end{theorem}

\begin{remark}\label{rem:ip->fv_starshaped}
Theorem \ref{thm:fv<->ip} opens the possibility to obtain upper estimates for the improved Poincar{\'e} constant $P_{\Omega}$ using known exact values or estimates for the corresponding Friedrichs-Velte and the Hardy constants.
Such upper estimates for $\Gamma_{\Omega}$ of a star-shaped planar or spatial domain are given in \cite{costabeldauge2015, horganpayne1983, payne2007} which upper estimates depend on the eccentricity of the domain with respect to the center of the star-shapedness.
As shown in \cite{chuawheeden2010} the Poincar{\'e} constant \eqref{ineq:iP} of a convex domain can be estimated by its eccentricity $\eta$, i.e. $P_{\Omega}\le c\eta^2$ for some positive constant $c$.
According to Theorem 6.2 in \cite{costabeldauge2015} and Theorem \ref{thm:fv<->ip} this remains valid for a planar star-shaped domain as well because its Hardy constant is at most 16, see \cite{ancona1986}.
\begin{equation}
P_{\Omega}\le 16\left(\eta+\sqrt{\eta^2-1}\right)^2\le 64\eta^2,
\end{equation}
where $\eta=\frac{R}{r}$ for the domain $\Omega$ star-shaped with respect to a disc of radius $r$ and contained in a concentric disc of radius $R$.
\hfill$\Box$
\end{remark}

As one sees the improved Poincar{\'e} inequalities imply the corresponding Friedrichs-Velte and Babu\v{s}ka-Aziz inequalities without any condition on the domain, the validity of the Hardy inequality is required only for the reverse implication.
Finiteness of the Hardy constant $H_{\Omega}$ was proved in \cite{necas1962} and \cite{boas1988} for Lipschitz and for H{\"o}lder domains, respectively.
However, as shown by Lehrb{\"a}ck \cite{lehrback2014}, not the smoothness of the boundary of the domain is decisive for $H_{\Omega}<\infty$ but rather its thickness in the sense of an inner boundary density condition formulated in Theorem 1.2 in \cite{lehrback2014}.
This reads in our special case \eqref{ineq:hardy} that we have $H_{\Omega}<\infty$ for an open set $\Omega\subset\mathbb{R}^n$ if the inner boundary density condition
\begin{equation}\label{ineq:contentcondition}
\mathcal{H}_{\infty}^{\lambda}\left(\partial\Omega\cap B(x,2d_{\Omega}(x))\right)\ge C_0 d_{\Omega}(x)^{\lambda}
\end{equation}
is satisfied for every $x\in\Omega$ with some positive constant $C_0$ and with an exponent $n-2<\lambda\le n-1$.
In \eqref{ineq:contentcondition} $B(x,r)$ denotes a ball centered in $x$ with radius $r$ and $\mathcal{H}_{\infty}^{\lambda}$ denotes the $\lambda$-dimensional Hausdorff content.

In particular each simply connected planar domain supports \eqref{ineq:hardy}, moreover, by \cite{ancona1986} we have $H_{\Omega}<16$ uniformly in the class of simply connected planar domains.
This implies that the Friedrichs-Velte, the Babu\v{s}ka-Aziz and the improved Poincar{\'e} inequalities are equivalent and we have
\begin{equation}\label{ineq:frveip_planesimplyconnected}
\frac{1}{4}\Gamma_{\Omega}\le P_{\Omega}\le 16\left(1+\Gamma_{\Omega}\right)
\end{equation}
for each simply connected planar domain, moreover, by \cite{jiang2014} they are equivalent to $\Omega$ being a John domain.
Altough this result is only a special case of that contained in \cite{jiang2014} but now we also have an explicit relation between the involved domain specific constants.
As shown in \cite{avkhadiev2014, barbatisetal2004, davies1995, davies1999} there are certain families of plane and spatial domains, the convex domains amongst them,  which have $H_{\Omega}=4$.
Using this along with \eqref{ineq:fviphconstants} there follows the more strict inclusion
\begin{equation}\label{ineq:frveip_convex}
\frac{1}{4}\Gamma_{\Omega}\le P_{\Omega}\le 4\left(1+\Gamma_{\Omega}\right),
\end{equation}
for these classes of domains, moreover, $\Gamma_{\Omega}\le P_{\Omega}\le 4\left(1+\Gamma_{\Omega}\right)$ for convex polygons by Remark \ref{rem:polygons}.

For general plane domains with more than one boundary points Avkhadiev \cite{avkhadiev2006} derived the bounds
\begin{equation}\label{ineq:avkhadievhardy}
M_0(\Omega)\le H_{\Omega}\le 4\left(\pi M_0(\Omega)+\frac{\Gamma(\frac{1}{4})^{\frac{1}{4}}}{4\pi^2}\right)^2,
\end{equation} 
where the quantity
\begin{equation}\label{eq:modulus}
M_0(\Omega)=\sup\frac{1}{2\pi}\log\frac{R(A)}{r(A)}
\end{equation}
is the maximum modulus of annuli $A$ centered on the boundary and separating the domain, i.e. $A=\left\{z\in\mathbb{C}:r(A)<|z-z_0|<R(A)\right\}\subset\Omega$ and $z_0\in\partial\Omega$.
That is, $H_{\Omega}<\infty$ iff $M_0(\Omega)<\infty$, which occurs if and only if $\partial\Omega$ is uniformly perfect.
Hence, in the family of plane domains with uniformly perfect boundary the Friedrichs-Velte, the Babu\v{s}ka-Aziz and the improved Poincar{\'e} constants are simultaneously finite or infinite and we have
\begin{equation}\label{ineq:frveip_planeuniformlyperfect}
\frac{1}{4}\Gamma_{\Omega}\le P_{\Omega}\le 
4\left(\pi M_0(\Omega)+\frac{\Gamma(\frac{1}{4})^{\frac{1}{4}}}{4\pi^2}\right)^2
\left(1+\Gamma_{\Omega}\right).
\end{equation}
This generalizes \eqref{ineq:frveip_planesimplyconnected} because we have $M_0(\Omega)=0$ for simply connected plane domains, however, the constant 16 in \eqref{ineq:frveip_planesimplyconnected} is less than the corresponding constant in \eqref{ineq:frveip_planeuniformlyperfect} for $M_0(\Omega)=0$. 

In case of arbitrary dimensional domains we have by \cite{mazya1985} that the Hardy constant $H_{\Omega}$ in \eqref{ineq:hardy} is finite if and only if there exists a positive constant $c$ such that
\begin{equation}\label{ineq:mazya}
\int_K\frac{1}{d_{\Omega}^2}\le c\cdot\text{\rm cap}_2(K,\Omega),
\end{equation}
for every compact subset $K\subset\Omega$, where the capacity of $K$ with respect to $\Omega$ is defined by
\begin{equation*}\label{eq:capacity}
\mathrm{cap}_2(K,\Omega)=\inf\left\{\int_{\Omega}|\nabla u|^2:
u\in C_0^{\infty}(\Omega)\text{ and }u(x)\ge 1\text{ for every }x\in K\right\},
\end{equation*}
moreover, we also have $c\le H_{\Omega}\le 16c$.
Substituting this into Theorem \ref{thm:fv<->ip} gives
\begin{equation}\label{ineq:capacity_estimate}
\frac{1}{4}\Gamma_{\Omega}\le P_{\Omega}\le
16\sup_{K\subset\Omega, K\mathrm{ compact}}\frac{\int_Kd_{\Omega}^{-2}}{\mathrm{cap}_2(K,\Omega)}
\left(1+\Gamma_{\Omega}\right)
\end{equation}
for the class of planar and spatial domains satisfying \eqref{ineq:mazya}.

If we want to obtain a more geometric condition for $\Omega$ than \eqref{ineq:mazya} in order to have $H_{\Omega}<\infty$, we can use the mean distance function $D_{\Omega}$ instead of the ordinary distance $d_{\Omega}$ as introduced by Davies \cite{davies1984}:
\begin{equation}\label{eqn:meandist}
\frac{1}{D_{\Omega}(x)^2}=\frac{1}{|\mathbb{S}^{n-1}|}\int_{\mathbb{S}^{n-1}}\frac{1}{d_{\nu}(x)^2}\mathrm{d}\nu,
\end{equation}
where $\mathbb{S}^{n-1}$ denotes the unit sphere in $\mathbb{R}^n$ and $d_{\nu}(x)=\inf\left\{|t|:x+t\nu\not\in\Omega\right\}$ for $x\in\Omega$ and $\nu\in\mathbb{S}^{n-1}$.
According to Theorem 17 in \cite{davies1984} we have the Hardy inequality
\begin{equation}\label{ineq:hardymean}
\int_{\Omega}\frac{u^2}{D_{\Omega}^2}\le\frac{4}{n}\int_{\Omega}|\nabla u|^2\text{ for every }u\in C_0^{\infty}(\Omega)
\end{equation}
for any domain of $\mathbb{R}^n$, $n\ge 2$.
One has $d_{\Omega}\le D_{\Omega}$ trivially and if $\partial\Omega$ fulfills some geometric condition then $D_{\Omega}$ can be estimated from above by some constant multiple of $d_{\Omega}$, see \cite{davies1984, davies1999, tidblom2004}.
For example if $\partial\Omega$ satisfies the exterior $\theta$-cone condition, e.g. each $x\in\partial\Omega$ is the vertex of an infinite circular cone of semi-angle $\theta$ which lies entirely in $\mathbb{R}^n\setminus\Omega$, then
we obtain
\begin{equation}\label{ineq:relate_distances}
d_{\Omega}\le D_{\Omega}\le 2d_{\Omega}\omega^{-\frac{1}{2}}\left(\frac{\sin\theta}{2}\right),
\end{equation}
and there follows $H_{\Omega}\le\frac{16}{n \omega\left(\frac{\sin\theta}{2}\right)}$ for $n\ge 2$, see \cite{davies1984}.
Here we have 
\begin{equation}\label{ineq:thetaconditionangle}
\omega\left(\alpha\right)=\frac{\int_0^{\arcsin\alpha}\sin^{n-2}t\,\mathrm{d}t}{2\int_0^{\frac{\pi}{2}}\sin^{n-2}t\,\mathrm{d}t}=
\begin{cases}
\frac{1}{\pi}\arcsin\alpha&\text{ if }n=2,\\
\frac{1-\sqrt{1-\alpha^2}}{2}&\text{ if }n=3,
\end{cases}
\end{equation}
for $0<\alpha<1$.
By Theorem \ref{thm:fv<->ip} for domains satisfying the exterior $\theta$-cone condition we obtain the estimate
\begin{equation}\label{ineq:thetacone_estimate3D}
\frac{1}{4}\Gamma_{\Omega}\le P_{\Omega}\le
\frac{16}{3\left(\frac{1}{2}-\sqrt{\frac{1}{4}-\frac{\sin^2\theta}{16}}\right)}
\left(1+\Gamma_{\Omega}\right),
\end{equation}
in the three dimensional case and
\begin{equation}\label{ineq:thetacone_estimate2D}
\frac{1}{4}\Gamma_{\Omega}\le P_{\Omega}\le
\frac{8\pi}{\arcsin\left(\frac{\sin\theta}{2}\right)}
\left(1+\Gamma_{\Omega}\right),
\end{equation}
in two dimensions.
These estimatates are more geometrical than \eqref{ineq:capacity_estimate}, however, for simply connected planar domains \eqref{ineq:frveip_planesimplyconnected} gives a better estimate then \eqref{ineq:thetacone_estimate2D}.

\subsubsection*{Concluding remarks}
The main result of this paper is twofold.
First, we proved by esimating the corresponding domain specific constants by each other that planar and spatial domains satisfying the Hardy inequality simultaneously support the Friedrichs-Velte inequality, the Babu\v{s}ka-Aziz inequality for the divergence and the improved Poincar{\'e} inequality for the gradient.
In the three dimensional case we derived with the same method a novel improved Poincar{\'e} inequality for the rotation.
As the geometrical considerations at the end of the paper show, the inequalities between the involved constants depend on the value of the Hardy constant of the domain for which one has many estimations available in the case of planar and spatial domains as well.

\end{document}